\documentclass{article}
\usepackage{geometry}                
\usepackage{graphicx}
\usepackage{tikz}
\usetikzlibrary{matrix,arrows,decorations.pathmorphing}
\usepackage{amssymb}
\usepackage{epstopdf}
\begin{document}
\newcommand*{\DashedArrow}[1][]{\mathbin{\tikz [baseline=-0.25ex,-latex, dashed,#1] \draw [#1] (0pt,0.5ex) -- (1.3em,0.5ex);}}%


\renewcommand{\thefootnote}{\fnsymbol{footnote}}
\newcommand{\starttext}{ \setcounter{footnote}{0}
\renewcommand{\thefootnote}{\arabic{footnote}}}
\renewcommand{\theequation}{\thesection.\arabic{equation}}
\newcommand{\be}{\begin{equation}}
\newcommand{\bea}{\begin{eqnarray}}
\newcommand{\eea}{\end{eqnarray}} \newcommand{\ee}{\end{equation}}
\newcommand{\N}{{\cal N}} \newcommand{\<}{\langle}
\renewcommand{\>}{\rangle}
\def\ba{\begin{eqnarray}}
\def\ea{\end{eqnarray}}
\newcommand{\PSbox}[3]{\mbox{\rule{0in}{#3}\includegraphics{#1}
\hspace{#2}}}

\def\v{\vskip .1in}

\def\al{\alpha}
\def\b{\beta}
\def\c{\chi}
\def\d{\delta}
\def\e{\epsilon}
\def\g{\gamma}
\def\l{\lambda}
\def\m{\mu}
\def\n{\nu}
\def\o{\omega}
\def\f{\phi}
\def\r{\rho}
\def\si{\sigma}
\def\t{\theta}
\def\z{\zeta}

\def\G{\Gamma}
\def\D{\Delta}
\def\O{\Omega}
\def\T{\Theta}

\def\cA{{\mathcal A}}
\def\cB{{\mathcal B}}
\def\cC{{\mathcal C}}
\def\cD{{\mathcal D}}
\def\cE{{\mathcal E}}
\def\cF{{\mathcal F}}
\def\cG{{\mathcal G}}
\def\cH{{\mathcal{H}}}
\def\cI{{\mathcal I}}
\def\cK{{\mathcal K}}
\def\cL{{\mathcal L}}
\def\cM{{\mathcal M}}
\def\cN{{\mathcal N}}
\def\cO{{\mathcal O}}
\def\cP{{\mathcal P}}
\def\cR{{\mathcal R}}
\def\cS{{\mathcal S}}
\def\cT{{\mathcal T}}
\def\cX{{\mathcal X}}

\def\N{\mathbb N}
\def\Z{{\mathbb Z}}
\def\Q{{\mathbb Q}}
\def\R{{\mathbb R}}
\def\C{{\mathbb C}}
\def\P{{\mathbb P}}
\def\K{{\rm K\"ahler }}

\def\KE{{\rm K\"ahler-Einstein }}
\def\KEE{{\rm K\"ahler-Einstein}}
\def\Rm{{\rm Rm}}
\def\Ric{{\rm Ric}}
\def\Hom{{\rm Hom}}
\def\mod{{\ \rm mod\ }}
\def\Aut{{\rm Aut}}
\def\End{{\rm End}}
\def\osc{{\rm osc}}
\def\vol{{\rm vol}}
\def\Vol{{\rm Vol}}
\def\ti\tilde

\def\u{\underline}

\def\pl{\partial}
\def\na{\nabla}
\def\i{\infty}
\def\I{\int}
\def\p{\prod}
\def\s{\sum}
\def\dd{{\bf d}}
\def\ddb{\partial\bar\partial}
\def\sub{\subseteq}
\def\ra{\rightarrow}
\def\hra{\hookrightarrow}
\def\Lra{\Longrightarrow}
\def\lra{\longrightarrow}
\def\LA{\langle}
\def\RA{\rangle}
\def\L{\Lambda}

\def\us{{\underline s}}
\def\Re{{\rm Re}}
\def\Im{{\rm Im}}
\def\tr{{\rm tr}}
\def\det{{\rm det}}
\def\half{ {1\over 2}}
\def\third{{1 \over 3}}
\def\ti{\tilde}
\def\un{\underline}
\def\Tr{{\rm Tr}}

\def\pz{\partial _z}
\def\pv{\partial _v}
\def\pw{\partial _w}
\def\w{{\bf w}}
\def\x{{\bf x}}
\def\y{{\bf y}}
\def\z{{\bf z}}
\def\tet{\vartheta}
\def\dwplus{\D _+ ^\w}
\def\dxplus{\D _+ ^\x}
\def\dzplus{\D _+ ^\z}
\def\chiz{{\chi _{\bar z} ^+}}
\def\chiw{{\chi _{\bar w} ^+}}
\def\chiu{{\chi _{\bar u} ^+}}
\def\chiv{{\chi _{\bar v} ^+}}
\def\os{\omega ^*}
\def\ps{{p_*}}

\def\hO{\hat\Omega}
\def\ho{\hat\omega}
\def\o{\omega}

\def\[{{\bf [}}
\def\]{{\bf ]}}
\def\Rd{{\bf R}^d}
\def\Ci{{\bf C}^{\infty}}
\def\pl{\partial}
\newcommand{\dotcup}{\ensuremath{\mathaccent\cdot\cup}}

\parindent=0in

\newtheorem{theorem}{Theorem}
\newtheorem{corollary}{Corollary}
\newtheorem{lemma}{Lemma}
\newtheorem{definition}{Definition}
\newtheorem{proposition}{Proposition}

\setcounter{equation}{0}

\centerline{ \bf \Large Convergence of the conical Ricci flow on $S^2$ to a soliton
\footnote{Work supported in part by
National Science Foundation grants DMS-12-66033, DMS-0847524 and DMS-0905873 and a Collaboration Grants for Mathematicians from Simons Foundation.} }
\bigskip\bigskip

\centerline{D.H. Phong*, Jian Song$^\dagger$,  Jacob Sturm$^\ddagger$ and Xiaowei Wang $^\ddagger$}

\bigskip

\smallskip

\medskip

\begin{abstract}

{\small In our previous work \cite{PSSW}, we showed that the
Ricci flow on $S^2$ whose initial metric has conical singularities $\s_j{\b_j}[p_j]$ converges
 to a constant curvature metric with conic singularities (in the stable and semi-stable cases)
or to a gradient shrinking soliton with conical singularities (in the unstable case). The purpose of this note
is to show that in the unstable case, that is, the case where $\b_k>\b_k'=\s_{j<k}\b_j$, that the limiting metric is the unique shrinking soliton with cone singularity
$\b_k[p_\i]+\b_k'[q_\i]$. This verifies the prediction made in \cite{PSSW}.
}

\end{abstract}
\section{Introduction}
Let $g_{S^2}$ be the round metric on $S^2$ and $\o_{S^2}= \sqrt{-1} g_{z\bar z}dz\wedge d\bar z$ the \K form, so that $[\o_{S^2}] = 2[p]$ for any $p\in S^2$.
Let $p_1,...,p_k\in S^2$ be a finite collection of points and $\b_1\leq \b_2\leq\cdots\leq\b_k \in (0,1)$. 
Let
$\b = \s_{j=1}^k \b_j[p_j]$. 
\v
A smooth metric $g$ on $S^2\backslash \{p_1,...,p_k\}$ is a cone metric
on $(S^2,\b)$ if it can be written in the form 
$$ \o\ =  \ 
{e^{f}\over \p_j|\si_j|^{\b_j}_{\o_{S^2}}}\cdot \o_{S^2}
$$ 
for some bounded function $f$ on $S^2$ where $\sigma_j$ is a section of $K^{-1}_{S^2}$ such that $[\si_j]=2[p_j]$. 
\v
A constant curvature  metric on $(S^2,\b)$ 
is a  metric $\o_\phi=\o_{S^2}+i\ddb\phi$ with the property

\be\label{eq} {\o_\phi }\ = \ {e^{-\g \phi}\over \p_j|\si_j|^{\b_j}_{\o_{S^2}}}\cdot \o_{S^2}
\ee
\v
where $\g = 1-{1\over 2}\s \b_j$. 
\v
An alternative form of (\ref{eq}) is

\be\label{eq1} \Ric(\o_\phi)\ = \ \g\o_\phi\ + \ \s\b_j[p_j]\ .
\ee

The Ricci flow is given by

\be\label{rf1} e^{-\dot\phi}\o_\phi\ = \ {e^{-\g \phi}\over \p_j|\si_j|^{\b_j}_{\o_{S^2}}}\cdot {\o_{S^2}}\ \ , \ \ \phi(0)=\phi_0
\ee
where $\dot\phi=\pl_t\phi$.
An alternative form of (\ref{rf1}) is

\be\label{ric} \pl_t g \ = \ -\Ric(g)\ + \ \g g\ , \ g(0)=g_0\ \ \hbox{\rm a  cone metric on $(S^2,\b)$} \ .
\ee
Here, and in all that follows,
we assume $\g>0$. We shall also assume that $g_0$ is ``regular" in the sense
of \cite{PSSW}. This means that $g_0$ is smooth on $S^2\backslash\{p_1,...,p_k\}$ and that in a neighborhood of $p_j$ there is a holomorphic
coordinate $z$ such that
\be\label{initial} g_0\ = \ e^{u_j}{idz\wedge d\bar z\over |z|^{2\b_j}}
\ee
where
$$ u_0,\,\D_0u_0,\, \D_0(\D_0u_0)\ \in \ C^{2,\al}(S^2,\b)\cap W^{1,2}\ .
$$
Here $C^{2,\al}(S^2,\b)$ is the Yin-H\"older space defined in \cite{Y}. In particular, if $u_0$ is harmonic in a neighborhood of $p_j$, then $g_0$
is regular.
\v
Let
$ \b_k' \ = \ \s_{j<k}\b_j$ and let 
$$ \b_\i = \b_k[p_\i] + \b_k'[q_\i]
$$
where
$p_\i$ and $q_\i$ are the north and south pole respectively.
Then we say $\b$ is stable, semi-stable or unstable if $\b_k'>\b_k$, $\b_k'=\b_k$, or $\b_k'<\b_k$ respectively.
\v
When $\b$ is stable it is known, by the work of 
\cite{MRS}, 
that the Ricci flow 
converges to the unique constant scalar metric on $(S^2,\b)$. In
\cite{PSSW} we give a new proof of this result. We also show that
in the semi-stable case, the Ricci flow converges to the unique constant scalar
curvature metric on $(S^2,\b_\i)$.

\v
We now assume that $\b$ is unstable and we let
$g_{\i}$ the unique conic shrinking soliton on $(S^2,\b_\i)$.
This means that $g_{sol}$ (which is rotationally symmetric by uniqueness) satisfies the following equation on $S^2\backslash\{p_\i,q_\i\}$:
 \begin{equation}
 R(g_{sol}) = \g + \Delta_{g_{sol}} \theta_{sol},  ~  \nabla_{g_{sol}}^2 \theta_{sol} = \frac{1}{2} (\Delta_{g_{sol}} \theta_{sol} ) g_{sol}, ~\int_{S^2}  e^{\theta_{sol}}dg_{sol} =2
 \end{equation}
  for a unique $\theta_{sol}\in 
  C^0(S^2)\cap C^\i(
  S^2\backslash\{p_\i,q_\i\}
  )$. Here $R(g_{sol})$ is the scalar curvature of $g_{sol}$.
  \v
We wish to prove the following:

\begin{theorem} For any initial regular metric $g_0$ on $(S^2,\b)$  the Ricci
flow  converges to $g_\i$.
\end{theorem}
Remark: In \cite{PSSW} we proved that there is a partition 
$\{1,2,...,k\}=I\cup J$ into disjoint subsets with the following property.
The flow (\ref{ric}) converges to the unique K\"ahler-Ricci solition $g_{I,J}$ with
cone structure

$$\b_{I,J}\ = \   \big(\s_{i\in I}\b_i\big)p_i\ + \ \big(\s_{j\in J}\b_j\big)p_j \ .
$$
Thus the content of Theorem 1 is that $I=\{p_k\}$ and $J=\{p_1,...,p_{k-1}\}$. In particular, $p_k\ra p_\i$ and $p_1,...,p_{k-1} \ra q_\i$ as $t\ra\i$.
\section{The proof}
\setcounter{equation}{0}
Let $\b=\s_{j=1}^k\b_j[p_j]$ and $g$ a cone metric on  $(S^2,\beta)$. Let
$f\in C^0(S^2)\cap W^{1,2}(S^2)$. 
We  define the normalized $W$-functional for the pair $(g,f)$ by the same expression as in the smooth case,
\begin{equation}
\label{W}
W(g, f) = \int_{S^2\setminus \beta } ({1\over 2\g}(R + |\nabla f|^2) + f){ e^{-f} \over 4\pi\tau} dg\ .
\end{equation}
We also define
$$   \mu (g)\ = \ \inf_f\,  \{W (g, f)\,:\,  \int_{S^2} e^{-f} dg =2\,\} 
\ .$$

Let $\m_1=\max\{\m(g_{I,J})\}$ and 
$\m_2=\max\{\m(g_{I,J})\,:\, g_{I,J}\not=g_\i$\} where, as above,
$$ g_\i\ = \ g_{I,J}\ , \ \hbox{where 
$I=\{p_k\}$ and $J=\{p_1,...,p_{k-1}\}$}
$$
We shall need the following (Lemma 7.3) from \cite{PSSW}
which was proved by first showing $\m(g_t)$ is increasing along the Ricci flow, and then using the toric structure of $g_{I,J}$ to compare $\m$ invariants.
\begin{lemma} We have $\m_1>\m_2$. Moreover, if there exists a regular cone metric $\ti g_0$ on 
$(S^2,\b)$ with
the property $\m(\ti g_0)>\m_2$, then the Ricci flow on $(S^2,\b)$ converges
in the Gromov-Haudsdorff $C^\i$ topology to $g_\i$ for any initial metric $g_0$. Thus $(S^2,g_t)\ra (S^2,g_\i)$ as metric spaces in the Gromov-Hausdorff topology. Moreover, for any compact subset $K\sub S^2\backslash\{p_\i,q_\i\}$ there exists a family of diffeomorphisms $f_t: S^2\ra S^2$ such that $f_t^*g_t\ra g_\i$ in $C^\i(K)$.
\end{lemma}

To prove Theorem 1,
we start by choosing coordinates on $\P^1$ in such a way that $p_\i$
is the point at infinity and $q_\i$ is the origin in $\C$. We define $g_\b$ to be the conic metric on
$(S^2,\b)$ whose \K form is given by

\be \o_\b\ = \ c(\b){\chi(z) dz\wedge d\bar z \over
\p_{j=1}^{k-1} |z-p_j|^{2\b_j}} + 
c(\b){(1-\chi(z))dz\wedge d\bar z\over (1+|z|^2)^{2-\b_k}
}
\ = \ 
F_\b\,\o_{FS}\ .
\ee

Here $\chi$ is smooth with compact support on $\C$ and equal to one in a large ball $B$ centered at $0\in\C$ which contains $p_1,...,p_{k-1}$  zero on the ball $2B$. The constant $c(\b)$ is chosen so that $\I dg_\b=2$ 
Thus $q_\i=0\in\C$ and $p_\i=\i\in\P^1$.
\v
 We have
$$  \Ric\, g_\b\ = \ \s_{j=1}^{k-1}\b_j[p_j]\ \  \hbox{on the
ball $B$}\ .
$$

Note that $c(\b)$ is a continuous function of $p_1,...,p_{k-1}$ and is thus bounded from above and away from zero provided $p_1,...,p_{k-1}$ remain in a bounded
subset of $\C$.
\v
Let $\r(t):\P^1\ra \P^1$ be the map $z\mapsto tz$ for $0\leq t\leq 1$ and
$\r(p_j)=p_j(t)$. Thus $p_j(t)=tp_j$ for $j<k$ and $p_k(t)=p_\i$.
Let $\b(t)=\s_{j=1}^k\b_j[p_j(t)]$. Then there exists $q>1$ such that 

\be\label{q}
\b(t)\ra \b_\i \hbox{\  and
$g_{\b(t)}\ra g_{\b_\i}$ in the GH topology and  
$F_{\b(t)}\ra F_{\b_\i}$ in $L^q$}\ .
\ee
Here and
in the following,
when we use notation such as $W^{1,2}, L^p,\D$ etc., the background metric is
always $g_{FS}$ unless otherwise specified.
\v
We wish to
construct a family of conic metrics $g_t$ on $(S^2,\b(t))$ such that
\be\label{conv} g_t\ra g_\i\ee
in the Gromov-Hausdorff topology 
 and
$\m(g_t)\ra \m(g_\i)$. 
Since $\m(g_\i)=\m_1>\m_2$ we conclude that for $t$ sufficiently large, 
$\m(g(t)) > \m_2$. For such a $t$, we let $g_0 = \r(t)^*g(t)$. Then
$g_0$ is a conic metric on $(S^2,\b)$ with the property $\m(g_0)>\m_2$,
and so Lemma 1 applies to give the desired conclusion.
\v
To define $g(t)$ we first write
$$  g_\i\ = \ e^{u_\i} g_{\b_\i}
$$
for some continuous function $u_\i$ which is smooth on $S^2\backslash\b_\i$. Theorem 1.1 of Datar-Guo-Song-Wang \cite{DGSW} shows $u_\i$ is a ``smooth $S^1$ invariant conic metric". This implies the $u_\i$ is smooth on $\C\backslash\{0\}$ and that there is a smooth
 $S^1$ invariant function $\ti u_\i$  on $B$ with the property
$u_\i(z)= \ti u_\i(w)$ where $|w|^2=|z|^{2(1-\b_k')}$. In particular, $u(z)$ has a Taylor expansion of the form

$$ u_\i(z) = a_0 + a_1|z|^{2(1-\b_k')}+ a_2|z|^{4(1-\b_k')}+ \cdots 
a_m|z|^{2m(1-\b_k')} + O(|z|^{2(m+1)(1-\b_k')})\ .
$$
In particular, there exist $C>0$ such that on $B$
\be\label{24}|u_\i(z)-a_0|\ + \  |z\pl_zu_\i|+ |z^2\pl_z\bar\pl_z u_\i|\ \ \leq \ C|z|^{2-2\b_k'}\ .
\ee
\v
We would like to define $g(t):=e^{u_\i}g_{\b(t)}$. This would satisfy (\ref{conv}) but doesn't quite work since $u_\i$ is not $C^2$ on the complement of 
$\b(t)$ so $e^{u_\i}g_{\b(t)}$ is not a regular metric in the sense of \cite{PSSW}. Instead we proceed as follows. Let $\psi=1-\chi$  which is
zero on $B$ and $1$ outside $2B$.  Define 
$$u_\i(t,z)= a_0 + \psi(z/t)(u_\i(z)-a_0) \ \ \hbox{if $t>0$}\ .
$$ 

Thus for each $t$ we see $u_\i(t,z)\in C^\i(S^2\backslash \{p_\i\})$ and $u_\i(t,z)$ is constant on the ball $tB$ and hence constant in a neighborhood of $p_1(t),....,p_{k-1}(t)$. Also, 
$$ u_\i(t,z) \ \longrightarrow\  u_\i(z)\ \ \ \hbox{pointwise as $t\ra 0$}
\ ,
$$

$$u_\i(t)\ra u_\i \ \hbox {and\ 
$\D_{g_{FS}}u_\i(t)\ra \D_{g_{FS}}u_\i$ uniformly on  compact subsets of $S^2\backslash \b_\i$}\ .
$$

Define 
$$ g(t) \ = \ e^{u_\i(t)} g_{\b(t)}\ .
$$
\v
We see that for each $t>0$ that $g(t)$ satsifies (\ref{initial}) with $u_j$ harmonic in a neighborhood of $p_j$. In particular, $g(t)$ is a regular metric. Moreover,
$$ \pl_z\bar\pl_z\, u_\i(t,z)\ = \ \psi'(z/t){1\over t}[\pl_zu_\i + \bar\pl_z u_\i]
\ +  \ 
\psi''(z/t){1\over t^2}(u_\i(z)-a_0))
\ + \ 
\psi(z/t)\pl_z\bar\pl_zu_\i\ .
$$
Since $|t|\geq c|z|$ when the right side is non-zero, we conclude from
(\ref{24}) that on $B$

$$ |\pl_z\bar\pl_z\, u_\i(t,z)|\ \leq \ {C\over |z|^{2\b_k'}}
$$
for some $C>0$ which is independent of $t$. We conclude that there exists
$q>1$ such that
$\|R(g(t))\|_{L^q}\leq C$ for all $t>0$.
Moreover, decreasing $q$ slightly if necessary,
 $\|R(t)\|_{L^q} \ra \|R(g_\i)\|_{L^q}$. 
 
 \v

In general, if $(X_t,g_t)\ra (X,g_\i)$ 
is any Gromov-Hausdorff limit of smooth manifolds, we know (\cite{Chow}, Lemma 6.28) that
$$ \m(g_\i) \geq \lim\sup \m(g_t)
$$
Thus our goal is to show
$ \lim_{T\ra\i}\inf_{t\geq T} \m(g(t))\ \geq  \ \m(g_\i)
$
 Assume not.
Then there exist $\d>0$ and a sequence $t_j\ra\i$ such that

\be\label{contra} \m(g_j)\ = \  \m(g(t_j))\ \leq \ \m(g_\i)-\d
\ee
Thus for each $j$ there is a positive function $\Phi_j=e^{-f_j/2}$ such that
$\|\Phi_j\|_{L^2(g_j)}=1$ and

\be\label{z} W(g_j, f_j)\ = \ \I_{S^2\backslash \b(t_j)} 
\big({2\over \g}|\na_j\Phi_j|^2 - \Phi_j^2{R_j\over 2\g} - 
\Phi_j^2\log\Phi_j^2\big)\, dg_j\ = \ \m(g_j)
\ \leq  \ 
\m(g_\i)-\d
\ee
\v
where $\g=1-{1\over 2}\s_j\b_j$.
\v

\begin{lemma} We have the following bounds.
\begin{enumerate} 
\item The $\Phi_j$ are uniformly bounded in $W^{1,2}$ that is, 
there exists $C>0$ such that
$$ \I \Phi_j^2\,dg_{FS}\ + \ \I \pl\Phi_j\wedge\bar\pl\Phi_j\ \leq \ C  \ \ \hbox{for all $j$}
$$ 
\item There exists $q>1$ such that  $\D u_j\ra \D u_\i$ in
$L^q$, that is
$$ \lim_{j\to\i} \I_{S^2} \, \big|\D_{g_{FS}}u_j-\D_{g_{FS}} u_\i)\big|^q\, dg_{FS}\ = \ 0
$$

\end{enumerate}
\end{lemma}
We postpone the proof  for the moment  and show how the lemma leads to a contradiction. 
\v
Part (1) implies
there exists $\Phi_\i \in W^{1,2}$ such that $\Phi_j \rightharpoonup \Phi_\i$ that is, $\Phi_j$ converges weakly to $\Phi_\i$ in $W^{1,2}$. Since
$W^{1,2}\hookrightarrow L^p$ is a compact imbedding for all $p>1$,
we see that after passing to a subsequence, $\Phi_j \ra \Phi_\i$ in $L^p$
for all $p>1$. Thus $\|\Phi_j\|_{L^p}\leq C_p$ for all $j$ and 
$\Phi_j^2\ra\Phi_\i^2$ in $L^p$ for all $p$.
\v
We claim

\be\label{1}
 \I \Phi_j^2{R_j} dg_j
\ = \ 
 \I \Phi_j^2(\D u_j + \g)\,dg_{FS}
 \ra \ 
 \I \Phi_\i^2(\D u_\i + \g)\,dg_{FS} 
 \ = \ 
\I \Phi_\i^2{R_\i} dg_j 
\ee

\be\label{11}1\ = \  \I \Phi_j^2 dg_j
\ = \
\I \Phi_j^2\,F_{\b_j}\,dg_{FS}
\ \ra \ 
\I \Phi_\i^2\,F_{\b_\i}\,dg_{FS}
\ = \ 
\I \Phi_\i^2 dg_\i
\ = \ 1
\ee

\be\label{2} \I \Phi_j^2\log\Phi_j^2 dg_j 
 =  
\I \Phi_j^2\log\Phi_j^2 F_{\b_j}dg_{FS}
 \ra 
\I \Phi_\i^2\log\Phi_\i^2 F_{\b_\i}dg_{FS}
 =  
\I \Phi_\i^2\log\Phi_\i^2 dg_j
\ee

\be\label{3} \lim\inf_j \I |\na_j\Phi_j|^2\,dg_j\ \geq \ 
\I |\na_\i\Phi_\i|^,dg_\i
\ee
To prove (\ref{1}) we note that $\D u_j \ra \D u_\i$ in $L^q$ and
$\Phi_j^2\ra \Phi_\i^2$ in $L^p$ for all $p$. Similarly (\ref{11}) 
follows from the fact that $F_{\b_j}\ra F_{\b_\i}$ in $L^q$ for some $q=q(\b)$.
\v

To prove (\ref{2}) we need only show $\Phi_j^2\log\Phi_j^2 \ra \Phi_\i^2\log\Phi_i^2 $ in $L^p$ for all $p$. To see this, first note that if $x,y>0$ there exists $\t$ between $x$ and $y$ such that

\be\label{mv} |x^2\log x^2- y^2\log y^2|\ = \ |4\t\log\t+2\t|\cdot |x-y|\ \leq \ 
C_\d (1+|x|^{2}+|y|^{2})\cdot |x-y|
\ee
by the mean value theorem (c.f. \cite{R}). Now substitute $x=\Phi_j$ and $y=\Phi_\i$ and apply
H\"older's inequality.

\v
Finally  (\ref{3}), which is equivalent to 
\be\label{6} \lim\inf_j \I\pl\Phi_j\wedge\bar\pl\Phi_j\ \geq \ 
\I\pl\Phi_\i\wedge\bar\pl\Phi_\i\\ ,
\ee
Since $\Phi_j \ \rightharpoonup \Phi_\i$ in $W^{1,2}$ we know
\be\label{7} \lim\inf_j \|\Phi_j\|_{W^{1,2}} \geq \|\Phi_\i\|_{W^{1,2}}
\ee
But $\Phi_j \ra \Phi_\i$ strongly in $L^2(g_{FS})$. Thus (\ref{6}) follows from 
(\ref{7}). 
\v
Taking the $\lim_{j\ra\i}$ of both sides of (\ref{z}) and applying
(\ref{1}),(\ref{11}),(\ref{2}) and (\ref{3}),  we obtain

\be\label{limj}  \ \I [{2\over \g}|\na_j\Phi_\i|^2 - \Phi_\i^2{R_\i\over 2\g} - \Phi_\i^2\log\Phi_\i^2] e^{u(t_\i)}g_{\b_\i}\ \leq  \   \ \m(g_\i)-\d
\ee
which contradicts the definition of $\m_(g_\i)$.
\v
Thus we have reduced the proof of Proposition 1 to the proof of the lemma.
\v
To prove the lemma, note that $(\ref{z})$ implies

$$ \I \pl\Phi_j\wedge\bar\pl\Phi_j\ \ \leq \ C\|\Phi_j\|^2_{L^q}
$$
for some $q>2$. On the other hand

$$ \|\Phi_j\|^2_{L^{q.}} - C_{q'}\|\Phi_j\|^2_{L^{2}}\leq \ C_{q'}\I \pl\Phi_j\wedge\bar\pl\Phi_j\
$$
for any $q'>q$.
This implies $\|\Phi_j\|_{L^q}\leq C$ and hence $\I \pl\Phi_j\wedge\bar\pl\Phi_j\leq C$. 
\v
This concludes the proof first part of the lemma. The second follows from the fact that $\D u_j \ra \D u_\i$
pointwise almost everywhere and $\|\D u_j\|_{L^q}$ is uniformly bounded for
some $q=q(\b)>1$.


\begin{thebibliography}{99}

{\small


\bibitem[Chow]{Chow} B. Chow et.al, ``The Ricci Flow: Techniques and Applications, Part I", Math. Surveys and Monongraphs, AMS (2007)

\bibitem[DGSW]{DGSW} V. Datar, B. Guo, J. Song and X. Wang,
``Connecting toric manifolds by conical Kahler-Einstein metrics",
arXiv:1308.6781

\bibitem[MRS]{MRS} Mazzeo, R., Y. Rubinstein and N. Sesum,
``Ricci flow on surfaces with conic singularities", 
arXiv:1306.6688


\bibitem[PSSW]{PSSW} Phong, D.H., J. Sturm, J. Song and X. Wang,
 ``The Ricci flow on the sphere with marked points", arXiv:1407.1118



\bibitem[R]{R} Rothaus, O.S.,
``Logarithmic Sobolev Inequalities and the Spectrum of Schr\"odinger Operators",
J. of Fun. Anal. {\bf 43} 110--120 (1981)
}


\bibitem[Y]{Y} Yin, Hao,
``Ricci flow on surfaces with conical singularities"  
J. Geom. Anal. 20 (2010), no. 4, 970--995. 

\end{thebibliography}
\end{document}